\documentclass[12pt]{amsart}

\usepackage{upref}
\usepackage{enumitem,dsfont,amssymb,bm}
\usepackage{amsmath}
\usepackage{amsthm}
\usepackage{graphicx,color}
\usepackage{mathtools}
\usepackage{mathrsfs}
\usepackage{braket}
\usepackage[english]{babel}
\usepackage[colorlinks=true, pdftitle={TITLE}, pdfauthor={Timothy 
	Allen Collier and Daniel Hauer}]{hyperref}

\numberwithin{equation}{section}

\newtheorem{theorem}{Theorem}[section]
\newtheorem{proposition}[theorem]{Proposition}

\newtheorem{corollary}[theorem]{Corollary}

\newtheorem{remark}[theorem]{Remark}

\newcommand\R{{\mathbb{R}}}

\newcommand\E{\mathcal{E}}

\newcommand\td{\mathrm{d}}
\newcommand\dx{\mathrm{d}x }

\DeclareMathOperator*{\divergence}{div}

\DeclareMathOperator{\loc}{loc}

\DeclareMathOperator{\Rg}{Rg}

\def\1{\raisebox{2pt}{\rm{$\chi$}}}
\DeclareMathOperator{\sign}{sign}

\newcommand\abs[1]{\lvert#1\rvert}
\newcommand\norm[1]{\lVert#1\rVert}

\definecolor{darkred}{rgb}{0.7,0.1,0.1}

\begin{document}

%\begin{frontmatter}
    \title[A density result for operators in $L^1$]{Density of the domain\\ of doubly nonlinear operators in $L^1$}
        
    \author{Timothy A. Collier}
    \author{Daniel Hauer}
    \address{School of Mathematics and Statistics, The University of Sydney, NSW 2006, Australia}
    \email{\href{mailto:timothyc@maths.usyd.edu.au}
    {timothyc@maths.usyd.edu.au}}
    \email{\href{mailto:daniel.hauer@sydney.edu.au}{daniel.hauer@sydney.edu.au}}

    \begin{abstract}
      The aim of this paper is to provide sufficient conditions implying
      that the effective domain $D(A\phi)$ of an $m$-accretive operator $A\phi$
      in $L^1$ is dense in $L^1$. Here,
      $A\phi$ refers to the
      composition $A\circ \phi$ in
      $L^1$ of the part $A=(\partial\E)_{\vert L^{1\cap \infty}}$ in
      $L^{1\cap\infty}$ of the subgradient $\partial\E$ in $L^2$ of a
      convex, proper, lower semicontinuous functional $\E$ on $L^2$ and
      a continuous, strictly increasing function $\phi$ on the real line
      $\R$. To illustrate the role of the sufficient conditions, we
      apply our main result to the class of doubly nonlinear operators
      $A\phi$, where $A$ is a classical Leray-Lions operator.
    \end{abstract}
    


\maketitle
\tableofcontents

\section{Introduction and main result}
%The density of the domain of a given operator plays an important role in the existence theory of associated Cauchy problems. 
%In this article we obtain such a density result for composed operators $A\varphi$ where $A$ is the subdifferential of an energy functional and $\varphi \in C(\R)$ is strictly increasing.

Let $(\Sigma,\mathcal{B},\mu)$ be a $\sigma$-finite measure
space, and $L^2:=L^2(\Sigma,\mu)$ the standard Lebesgue space
of $\mu$-a.e. equivalence classes of square-integrable
functions $u : \Sigma\to \R$ endowed with the $L^2$-inner
product 
\begin{displaymath}
    \langle u,v\rangle_{L^2}:=
    \int_{\Sigma}u\,v\,\td\mu\qquad
    \text{for every }u,\,v\in L^2,
\end{displaymath}
and the induced norm $\norm{u}_2:=\sqrt{\langle
u,u\rangle_{L^2}}$. Further, let $\E : L^2 \to (-\infty,
+\infty]$ be a proper, lower semicontinuous, convex functional,
$D(\E):=\{u\in L^2\,\vert\, \E(u)<\infty\}$ the effective
domain of $\E$, and $A:=\partial \E$ the \emph{sub-gradient} in $L^2$
of $\E$ defined through the graph
\begin{displaymath}
    \partial\E=\Big\{(u,h)\in D(\E)\times L^2\,\Big\vert\,
    \langle h,v-u\rangle_{L^2}\le \E(v)-\E(u)\,\forall\,v\in
    L^2\Big\}.
\end{displaymath}
Then the operator $A$ is \emph{maximal monotone} (or, equivalently,
\emph{$m$-accretive}) in $L^2$, that is, for every $\lambda>0$, $A$
satisfies for $X=L^2$ the so-called \emph{range condition}
\begin{equation}
  \label{eq:range}
  \Rg\big(I_{L^2}+\lambda\,A\big)=X
\end{equation}
and the \emph{resolvent operator}
$J_{\lambda}:=(I_{L^2}+\lambda A)^{-1}$ is a contraction mapping on
$L^2$. In particular, $-A$ generates a $C_0$-semigroup
$\{e^{-tA}\}_{t\ge 0}$ of contractions on the closure $\overline{D(A)}$
in $L^2$ of the domain $D(A):=\{u\in L^2\,\vert\,Au\neq\emptyset\}$ of
$A$, and for $A$ the following density result holds (see, for instance,
\cite{MR2582280} or \cite{MR0348562}).

\begin{theorem}[{Density of $A$ in $L^2$}]
    \label{thm:classicalDensity}
    Let $\E : L^2 \to (-\infty,\infty]$ be a proper, lower
    semicontinuous, convex functional, and $A=\partial\E$ the
    sub-gradient in $L^2$ of $\E$. Then, the domain $D(A)$ is a dense
    subset of $D(\E)$ in $L^2$.
\end{theorem}

It is worth noting that the statement of
Theorem~\ref{thm:classicalDensity} remains true for sub-differential
operators $\partial\E\subseteq X\times X^{\ast}$ of a proper, lower
semicontinuous, convex functional $\E : X \to (-\infty,+\infty]$ defined
on a real Banach space $X$ (see, for instance, \cite[Proposition
1.6]{MR2582280}). But for our purposes here, it is sufficient to
consider $X=L^2$ or $L^1$ the standard Lebesgue space
$L^1(\Sigma,\mu)$.\medskip

Accordingly to the classical (nonlinear) semigroup theory in Hilbert
spaces (see, for instance, \cite{MR0348562}), for every initial value
$u_0\in \overline{D(A)}$, $u(t):=e^{-tA}u_0$, $t\ge0$, is a (strong)
solution of the Cauchy problem (in $L^2$)
\begin{equation}
	\label{eq:1}
	\begin{cases}
    \begin{alignedat}{2}
	\dot{u}(t)+Au(t)&\ni 0\quad && 
	\text{on $(0,\infty),$}\\
	u(0)&=u_{0}.\quad && 
    \end{alignedat}
	\end{cases}
\end{equation}
Hence for application, it is worth knowing which functions $u_0$
belong to the closure $\overline{D(A)}$ in $L^2$. But thanks to
the density result given by Theorem~\ref{thm:classicalDensity},
the characterization of the set $\overline{D(A)}$ is obtained
via the density of the effective domain $D(\E)$ of the
functional $\E$, which, in practice, is much easier to show
than to characterize directly the set $\overline{D(A)}$ (see,
for instance, \cite{MR2033382,MR4041276,MR3465809,CoulHau2016,MR3369257,MR4200826,Hauer2019,MR3491533,MR4365127,MR3427974,MR3057171,MR2289546}
and note that this is certainly not a complete list).\medskip

%The density of the domain of such an operator plays an
%important role in the existence theory of associated Cauchy
%problems.
%Theorem~\ref{thm:classicalDensity} plays an important role
%when analyzing equations driven by a subdifferential operator,
%such as the $p$-Laplacian or fractional $p$-Laplacian.
%It is standard to consider the well-posedness of associated
%abstract Cauchy problems of the form
%\begin{equation}
%	\label{eq:1}
%	\begin{cases}
%    \begin{alignedat}{2}
%	\dot{u}(t)+Au(t)&\ni g(t)\quad && 
%	\text{on $(0,T),$}\\
%	u(0)&=u_{0}\quad && \text{in $X$,}
%    \end{alignedat}
%	\end{cases}
%\end{equation}
%where $A:X\rightarrow X$ is a quasi-accretive operator on a
%Banach space $X$, the function $g \in L^1(0,T;X)$ and $u_0 \in
%X$. 
%If $A$ satisfies a further \emph{range condition}, then by the
%existence theorem of Crandall and Liggett \cite{MR0287357} we
%have existence of mild solutions for all $u_0 \in
%\overline{D(A)}^{\mbox{}_{X}}$. 
%Such solutions may then be regularized to find strong
%solutions.
%We refer to the monograph \cite{MR2582280} for the general
%theory. 
%Hence, in this setting, density of $D(A)$ in $X$ provides
% well-posedness in the sense of mild solutions for all initial
%data %$u_0 \in X$. 
%Note that if $X$ is reflexive then we can interpret $\partial
%\E$ as an operator from $X \rightarrow X$ to consider $A =
%\partial \E$, hence it is common to set $X$ to be the 
%Lebesgue space $L^2$.

The aim of this article is the present a generalization of this
density result (Theorem~\ref{thm:classicalDensity}) for the
perturbed operator $A\phi$ in $L^1$ defined by
\begin{displaymath}
    A\phi = \Big\{(u,v)\in L^{1}\times L^{1}\,\Big\vert\, 
    \phi(u) \in D(A) \text{ and } v \in A(\phi(u))
    \Big\},
\end{displaymath}
where the operator $A=(\partial\E)_{\vert L^{1\cap\infty}}$ the
\emph{part}/\emph{restriction} of the sub-gradient $\partial\E$
in $L^{1\cap\infty}:=L^1\cap L^{\infty}$ given by
\begin{displaymath}
    (\partial\E)_{\vert L^{1\cap\infty}}=\partial\E\cap 
    (L^{1\cap\infty}\times L^{1\cap \infty}),
\end{displaymath}
and the graph (in $L^1$)
\begin{displaymath}
    \phi:=\Big\{(u,h)\in L^1\times L^1\,\Big\vert\,
    h(x)=\phi(u(x))\quad\text{$\mu$-a.e. on }\Sigma\Big\}
\end{displaymath}
of a strictly increasing, surjective, and continuous function
$\phi : \R \to \R$ satisfying $\phi(0)=0$. The operator $A\phi$
is then nothing less than the composition $A\circ \phi$
realized in $L^1$ of the sub-gradient $\partial\E$ and $\phi$.
Moreover, under the additional assumptions that 
\begin{align}
    \label{eq:Hyp-A0}\tag{H1}
    &(0,0)\in A\text{ (or, equivalently, $\E(0)=\min_{L^2}
    \E$);}\\ \label{eq:Hyp-A-compl-acc}\tag{H2}
    &A \text{ is \emph{completely accretive}, that is, }\\
    \label{eq:integral-complete-accretive}
    &\qquad\int_{\Sigma}j(u-\hat{u})\,\td\mu
    \le \int_{\Sigma}j(u-\hat{u}+\lambda(v\hat{v}))\,\td\mu\\
    \notag
    &\text{for every $(u,v)$, $(\hat{u},\hat{v})\in A$, and
    $j\in \mathcal{J}_0$,}
    \intertext{where $\mathcal{J}_0 := \{j : \R\to [0,\infty]\,
    \vert\,j\text{ convex, lower semicontinuous, }j(0)=0\}$,} 
    \label{eq:sum-rule}\tag{H3}
    &\text{the Yosida operator $[\phi^{-1}]_{\lambda}$ of the 
    inverse graph $\phi^{-1}$ given by}\\ \notag
    &\qquad [\phi^{-1}]_{\lambda}=\lambda^{-1}
      \Big(I_{L^1}-(I_{L^1}+\lambda \phi^{-1})^{-1}\Big)\\ \notag
    &\text{satisfies}\\ \label{eq:H3-Yosida-1}
    &\qquad \int_{\Sigma} v\,[\phi^{-1}]_{\lambda}u\,\td\mu\ge 0
      \qquad\text{for every $(u,v)\in A$, and}\\ \label{eq:H3-Yosida-2}
    &\qquad \int_{\Sigma} v\,\sign_0\Big((\phi^{-1})_{\lambda}u\Big)\,
      \td\mu\ge 0\qquad\text{for every $(u,v)\in A$,}
\end{align}
where $\sign_0 : \R\to \{-1,0,1\}$ is the \emph{restricted
signum} function given by  
\begin{displaymath}
    \text{sign}_0(s) := \begin{cases}
        1& \text{if } s > 0,\\
        0& \text{if } s = 0,\\
        -1& \text{if } s < 0,
    \end{cases}\qquad\text{for every $s\in \R$,}
\end{displaymath}
Coulhon and Hauer~\cite[Proposition~2.18]{CoulHau2016} showed that the
closure $\overline{A\phi}$ of $A\phi$ in $L^1$ is
\emph{$m$-$T$-accretive} in $L^1$ with a \emph{complete resolvent}, that
is, $\overline{A\phi}$ satisfies for $X=L^1$ the range
condition~\eqref{eq:range}, and for every $\lambda>0$, the resolvent
$J_{\lambda}$ of $\overline{A\phi}$ is \emph{$T$-contractive} (on $L^1$)
\begin{displaymath}
  \norm{[J_{\lambda}u-J_{\lambda}u_2]^{+}}_1\le \norm{[u_1-u_2]^{+}}_1
\end{displaymath}
for every $u_1$, $u_2\in L^1$, and \emph{complete} in the sense
\begin{displaymath}
 \int_{\Sigma}j(J_{\lambda}u)\,\td\mu
    \le \int_{\Sigma}j(u)\,\td\mu
\end{displaymath}
for every $u_2\in L^{1\cap\infty}$ and $j\in \mathcal{J}_0$. In particular, $-\overline{A\phi}$
generates a $C_0$-semigroup $\{e^{-t \overline{A\phi}}\}_{t\ge 0}$ of
$T$-contractions\footnote{Or, equivalently, $\{e^{-t
    \overline{A\phi}}\}_{t\ge 0}$ is contractive on the closure
  $\overline{D(A\phi)}$ in $L^1$ and \emph{order-preserving}, that is,
  for every $u_1$, $u_2\in \overline{D(A\phi)}$ satisfying $u_1\le u_2$
  $\mu$-a.e. implies that for all $t\ge 0$, $e^{-t
    \overline{A\phi}}u_1\le e^{-t
    \overline{A\phi}}u_2$ $\mu$-a.e..} $e^{-t \overline{A\phi}}$ on the closure
$\overline{D(A\phi)}$ in $L^1$.  Thus, under the hypotheses
\eqref{eq:Hyp-A0}--\eqref{eq:sum-rule}, the Cauchy problem~\eqref{eq:1}
(in $L^1$) for $A$ replaced by $\overline{A\phi}$ is well-posed.\medskip

We emphasize that \cite[Proposition~2.18]{CoulHau2016} provides a
nonlinear realization of the construction of an $m$-$T$-accretive
operator $A\phi$ in $L^1$ by Crandall and Pierre~\cite{MR647071}, who
studied the case when $A$ is \emph{linear}. But they did not show that
the domain $D(A\phi)$ is dense in $L^1$. Now, under the hypotheses
\eqref{eq:Hyp-A0}--\eqref{eq:sum-rule} and if the effective domain
\begin{displaymath}
    D(\E\phi) := \Big\{u \in L^{1\cap\infty}\,\Big\vert\,
    \phi(u) \in D(\E)\Big\}
\end{displaymath}
satisfies the \emph{lattice property}
\begin{equation}
    \label{eq:latice}\tag{H4}
    \text{for every $u\in D(\E\phi)$, one that $u^{+}$ and
    $u^{-}\in D(\E\phi)$,}
\end{equation}
where $u^{+}:=\max\{u,0\}$ and $u^{-}:=\min\{u,0\}$ denote the
\emph{positive} and \emph{negative} part of a function 
$u\in L^{1}_{\loc}$, then our main results provides the
following density result for the perturbed operator $A\phi$
holds.

\begin{theorem}[{Density of $A\phi$ in $L^1$}]
    \label{thm:L1densityofE}
    Let $\E : L^2 \to (-\infty,\infty]$ be a proper, lower
    semicontinuous, convex functional, and $\phi\in C(\R)$ be a strictly
    monotone increasing, and surjective function satisfying
    $\phi(0)=0$. Suppose $A=(\partial\E)_{\vert L^{1\cap\infty}}$ and
    $\phi$ satisfy the hypotheses
    \eqref{eq:Hyp-A0}--\eqref{eq:latice}. Then the domain $D(A\phi)$ is
    a dense subset of $D(\E\phi)$ in $L^1$.
\end{theorem}

In specific cases of $A$ and $\phi$, the density of $D(A\phi)$ in $L^1$
has been proved: for example, Evans \cite[Sect. 2, Proposition
1]{MR454360} for $A=-\Delta$ in $L^1$ and $\phi: \R\to\R$ increasing and
$\phi^{-1}$ Lipschitz continuous, or Igbida \cite[Proposition
2.1]{alma991006084359706657} for
$A=-\divergence(\abs{\nabla \cdot}^{p-2}\nabla\cdot)$ in
$L^{1\cap\infty}$ and $\phi(r)=\abs{r}^{m-1}r$ for $r\in \R$ and
$m\ge 1$. We outline the proof of Theorem~\ref{thm:L1densityofE} in the
subsequent section, and apply this result in
Section~\ref{sec:applications} to $A=(-\Delta_p)^s$ the fractional
$p$-Laplacian in $L^{1\cap\infty}$ for $1\le p<\infty$ and $0<s<1$, and
$A=\Delta_1$ the $1$-Laplacian in $L^{1\cap\infty}$.\medskip

In \cite[Lemma~7.1]{BenilanCompletelyAccretiveMR1164641},
B\'enilan and Crandall showed that if the functional $\E$
satisfies 
\begin{equation}
    \label{eq:complete-accretive-phi}\tag{H2$\mbox{}^{\ast}$}
    \E(u+T(\hat{u}-u))+\E(\hat{u}-T(\hat{u}-u))\le
        \E(u)+\E(\hat{u}),
\end{equation}
for every $u$, $\hat{u}\in L^2$ and $T\in \mathcal{P}_0$, 
where $\mathcal{P}_0$ denotes the set of all $T\in C(\R)$
satisfying $T(0)=0$, $0\le T'\le 1$ and $T'$ has compact
support, then $A$ satisfies
hypothesis~\eqref{eq:Hyp-A-compl-acc}. In application, 
hypothesis~\eqref{eq:complete-accretive-phi} is easier to
verify than the integral
condition~\eqref{eq:integral-complete-accretive}. In fact, both
hypotheses~\eqref{eq:Hyp-A-compl-acc} and
\eqref{eq:complete-accretive-phi} are related to invariance
principles satisfied by the semigroup $\{e^{-tA}\}_{t\ge 0}$
generated by $\E$ (see, for example, \cite{MR3465809}). Now, we
can conclude from Theorem~\ref{thm:L1densityofE} the following
corollary.

\begin{corollary}
\label{cor:densityResult}
Let $\E : L^2 \to (-\infty,\infty]$ be a proper, lower semicontinuous,
convex functional, and $\phi\in C(\R)$ be a strictly monotone
increasing, and surjective function satisfying $\phi(0)=0$. Suppose
$A=(\partial\E)_{\vert L^{1\cap\infty}}$ and $\phi$ satisfy the
hypotheses \eqref{eq:Hyp-A0}, \eqref{eq:complete-accretive-phi},
\eqref{eq:sum-rule}, and \eqref{eq:latice}. Then the domain $D(A\phi)$
is a dense subset of $D(\E\phi)$ in $L^1$.
\end{corollary}

Now, we turn now to the proof of our main result.

\section{Proof of the main result}
\label{sec:proof-main-result}

In fact, our proof of Theorem~\ref{thm:L1densityofE} generalizes the
simple proof of \cite[Proposition 1.2 given on p.48]{MR2582280}. In
order to highlight the main idea of our proof, we intend to briefly
recall the proof of Theorem~\ref{thm:classicalDensity}.

\begin{proof}[Reminder of the proof of
Theorem~\ref{thm:classicalDensity}]
   Let $f\in D(\E)$. Since $A$ is maximal monotone, 
   for every $\lambda>0$, $u_{\lambda}:=J_{\lambda}f$ is the
   unique solution of
   \begin{equation}
    \label{eq:classic-inclusion}
       u_{\lambda}+\lambda\,A u_{\lambda}\ni f,
   \end{equation}
   meaning that 
   \begin{displaymath}
       u_{\lambda}\in D(A)\quad\text{ and }\quad
       f-u_{\lambda} \in \lambda\partial
   \E(u_{\lambda}).
   \end{displaymath}
   Thus, by applying the definition of the sub-gradient, one obtains that
   \begin{displaymath}
       \norm{u_{\lambda}-f}_{2}+\lambda\,\big(\E(u_{\lambda})
       -\E(f)\big)\le 0.
   \end{displaymath}
   Since $\E$ is bounded from below by an affine function
   (cf~\cite[Proposition 1.1]{MR2582280}), sending $\lambda\to 0+$ in
   the latter inequality yields that $u_{\lambda}\to f$ in $L^2$ as
   $\lambda\to 0+$.
\end{proof}

Next, we turn to the proof of Theorem~\ref{thm:L1densityofE}. With the
previous proof in our mind, the first step is to find for given
$f\in D(\E\phi)$ and $\lambda>0$, a solution $u_\lambda$ of the
inclusion~\eqref{eq:classic-inclusion} with $A$ replaced by the
perturbed operator $A\phi$. In fact, the \cite[proof of
Proposition~2.18, pp.113]{CoulHau2016} provides an existence result of the
inclusion~\eqref{eq:classic-inclusion} under slightly more general
hypotheses than \eqref{eq:Hyp-A0}--\eqref{eq:sum-rule}. More
precisely, we have the following.

\begin{proposition}[{\cite{CoulHau2016}}]\label{prop:existence-Aphi}
  Let $A$ be an $m$-completely accretive operator on $L^2$ satisfying
  $(0,0)\in A$, and $\phi\in C(\R)$ be a surjective, strictly increasing
  function satisfying $\phi(0)=0$. Further, suppose that
  $A_{\vert L^{1\cap\infty}}$ and $\phi$ satisfy
  \eqref{eq:sum-rule}. Then, for every $\lambda>0$,
  $f\in L^{1\cap\infty}$, and every sufficiently small $\nu>0$, there is
  a unique $u_{\lambda}\in D(A_{\vert L^{1\cap\infty}}\phi)$ satisfying
  \begin{equation}
    \label{eq:rangeProblemAphi}
    u_{\lambda}+\lambda\,\big(\nu\,\phi(u_{\lambda})+
    A_{\vert L^{1\cap\infty}}\phi(u_{\lambda})\big)\ni f.
  \end{equation}
\end{proposition}

Now, we can give the proof of Theorem~\ref{thm:L1densityofE}, where
we apply Proposition~\ref{prop:existence-Aphi} to $A=(\partial\E)_{\vert_{L^{1\cap\infty}}}$.

\begin{proof}[Proof of Theorem~\ref{thm:L1densityofE}]
  Let $f\in D(\E\phi)$, $\lambda>0$, and $\nu>0$ sufficiently small. Then, by
  Proposition~\ref{prop:existence-Aphi}, 
%  It follows from the definition that $D(\partial \E) \subseteq D(\E)$
%  and by our definition of composition in $L^1$, that
%  $D((\partial\E)_{\vert_{L^{1\cap\infty}}}\varphi) \subseteq
%  L^{1\cap\infty}$.
%  Hence
%  $D((\partial\E)_{\vert_{L^{1\cap\infty}}}\varphi) \subseteq
%  D(\E\varphi)\cap L^{1\cap\infty}$.
%  Under the hypotheses of this theorem, we can then apply \cite[Lemma
%  A.3.1]{CoulHau2016} to $(\partial\E)_{\vert_{L^{1\cap\infty}}}\varphi$
%  and obtain that for every $\lambda>0$, 
%  $u \in D(\E \varphi)\cap L^{1\cap\infty}$, 
  there is a unique $u_\lambda \in D(A\phi)$
  satisfying~\eqref{eq:rangeProblemAphi}, or equivalently,
    \begin{displaymath}
       % \label{eq:rangeProblemAphi}
     \phi(u_{\lambda})\in D(\E)\quad\text{ and }\quad	   
     f-u_{\lambda} -\lambda\,\nu\,\phi(u_\lambda) \in  
           (\partial\E)_{\vert_{L^{1\cap\infty}}}\Big(\phi(u_\lambda)\Big).
    \end{displaymath}
    Thus, by applying the definition of the sub-gradient, one obtains
    that
    \begin{displaymath}
      \langle f-u_{\lambda} -\lambda\,\nu\,\phi(u_{\lambda}),
      \phi(f)-\phi(u_{\lambda}) \rangle_{L^2}
      \le \lambda\Big(\E(\phi(f))-\E(\phi(u_{\lambda}))\Big), 
    \end{displaymath}
    or, equivalently,
    \begin{equation}
      \label{eq:u-lambda-f-integral}
    \begin{split}
      \langle f-u_{\lambda},\phi(f)-\phi(u_{\lambda}) \rangle_{L^2}
     &+\lambda\Big(\E(\phi(u_{\lambda}))-\E(\phi(f))\Big)\\
      &\qquad \le \lambda\,\nu\,\langle \phi(u_{\lambda}),
      \phi(f)-\phi(u_{\lambda}) \rangle_{L^2}.
    \end{split}
  \end{equation}
%    \begin{align*}
%          \big (u_\lambda-u,\varphi(u_\lambda)-\varphi(u) \big)_{L^{2}}
%          =
%          & -\lambda\varepsilon \big(\varphi(u_\lambda),\varphi(u_\lambda)-\varphi(u) \big)_{L^{2}}\\
%          &\qquad -\lambda\big(
%          v_\lambda,\varphi(u_\lambda)-\varphi(u)\big)_{L^{2}}.
%   \end{align*}
    Since $\nu\phi(\cdot)+A\phi$ has a complete resolvent, one has that
    \begin{equation}
      \label{eq:bound-u-lambda}
      \int_{\Sigma}j(u_{\lambda})\,\td\mu
      \le \int_{\Sigma}j(f)\,\td\mu
    \end{equation}
    for every $j\in \mathcal{J}_0$ and $\lambda>0$, implying that
    \begin{displaymath}
      \norm{u_{\lambda}}_{\infty}\le \norm{f}_{\infty}
    \end{displaymath}
    for every $\lambda>0$. Thus, the first term on
    the right-hand side can be estimated by
	\begin{align*}
         % \begin{split}
         % -\big(\varphi(u_\lambda),\varphi(u_\lambda)-\varphi(u)
         % \big)_{L^{2}}
          \langle \phi(u_{\lambda}),
      \phi(f)-\phi(u_{\lambda}) \rangle_{L^2}&= - \langle \phi(f)-\phi(u_{\lambda}),
      \phi(f)-\phi(u_{\lambda}) \rangle_{L^2}\\
       &\hspace{2.5cm}  + \langle \phi(f),
      \phi(f)-\phi(u_{\lambda}) \rangle_{L^2}\\
        &\le \norm{\phi(f)}_{2}^2-\langle \phi(f),\phi(u_{\lambda})
          \rangle_{L^2}\\
         &\le \norm{\phi(f)}_{2}^2+\norm{\phi(u_\lambda)}_\infty\norm{\phi(f)}_{1}\\   
     %  & \le \norm{\varphi(u_\lambda)}_\infty\norm{\varphi(u)}_{1}\\
          &\le \norm{\phi(f)}_{2}^2+
            \sup_{[-\norm{f}_\infty,\norm{f}_\infty]}\norm{\phi}_\infty\norm{\phi(f)}_{1}.
    %	\end{split}
	\end{align*}
        According to~\eqref{eq:bound-u-lambda},
        $(u_{\lambda})_{\lambda>0}$ is bounded in $L^2$, and $\E$ is
        bounded from below by an affine function. Further, since $\phi$
        is monotonically increasing,
        \begin{displaymath}
          	f_{\lambda}(x):=\big(f(x)-u_{\lambda}(x)\big)
                \big(\phi(f(x))-\phi(u_{\lambda}(x)\big)
                \ge 0
        \end{displaymath}
        $\mu$-a.e. on $\Sigma$. Therefore, sending $\lambda\to 0+$
        in~\eqref{eq:u-lambda-f-integral} yields that
         % $u_{\lambda}\to f$ in $L^2$ as
         % $\lambda\to 0+$.
         % and since $v_{\lambda}\in
         % (\partial\E)_{\vert_{L^{1\cap\infty}}}\varphi (u_\lambda)$ 
         % and noting that $\E(\varphi(u_\lambda)) \ge -M$, it follows that
         % \begin{displaymath}
         %   \begin{split}
         %     -\big(
         %     v_{\lambda},\varphi(u_\lambda)-\varphi(u)\big)_{L^{2}}& =
         %     -\big(
         %     (\partial\E)_{\vert_{L^{1\cap\infty}}}\varphi(u_\lambda),
         % \varphi(u_\lambda)-\varphi(u)\big)_{L^{2}}\\
         %     & \le
         %     -\left(\E(\varphi(u_\lambda))-\E(\varphi(u))\right)\\
         %     & \le \E(\varphi(u))+M
         %   \end{split}
         % \end{displaymath}
         % for all $\lambda > 0$.
         % Thus and since $\varphi$ is increasing, we have shown that
         % \begin{align*}
         %               %   0\le
         %   (u_{\lambda}-u,\varphi(u_{\lambda})-\varphi(u))_{L^{2}}&\le
         %   \lambda\varepsilon\sup_{[-\norm{u}_\infty,\norm{u}_\infty]}
         %   \norm{\varphi}_{\infty}\norm{\varphi(u)}_{1}\\
         %         %   &\quad +\lambda\, \left(\E(\varphi(u))+ M\right)
         % \end{align*}
         % for all $\lambda>0$, from where we can conclude that
	\begin{displaymath}
	\lim_{\lambda\rightarrow
          0^{+}}f_{\lambda} = 0\qquad\text{in $L^1$.}
	\end{displaymath}
        Now, let $(\lambda_{n})_{n\ge 1}$ be an arbitrary zero sequence
        with $\lambda_{n}>0$ for every integer $n\ge 1$. Then, there is
        a subsequence $(\lambda_{k_{n}})_{n\ge 1}$ of
        $(\lambda_{n})_{n\ge 1}$ such that
	\begin{displaymath}
	\lim_{n\to \infty}f_{\lambda_{k_n}}(x) = 0
	\end{displaymath}
        for $\mu$-a.e.~$x \in \Sigma$. Due to the strict
        monotonicity of $\phi$, the letter implies that
        \begin{displaymath}
          \lim_{n\to \infty}u_{\lambda_{k_{n}}}(x) = f(x)
        \end{displaymath}
        for $\mu$-a.e.~$x \in \Sigma$. Since $(\lambda_{n})_{n\ge 1}$
        was an arbitrary zero sequence, we have thereby shown that for
        every $f\in D(\E\phi)$, there is a sequence
        $(u_{\lambda})_{\lambda>0}\subseteq D(A\phi)$ such that 
	\begin{equation}
        \label{eq:pointwise-mu-conv}
	\lim_{\lambda\rightarrow 0^{+}}u_{\lambda}(x) = f(x) 
        \qquad\text{for $\mu$-a.e.~$x \in
          \Sigma$.}
	\end{equation}
        To conclude this proof, it remains to show that for every zero
        sequence $(\lambda_{n})_{n\ge 1}$ with $\lambda_{n}>0$, there is
        a subsequence $(\lambda_{k_{n}})_{n\ge 1}$ of
        $(\lambda_{n})_{n\ge 1}$ and $g\in L^1$ satisfying
        \begin{displaymath}
          \abs{u_{\lambda_{k_{n}}}}\le g\qquad\text{$\mu$-a.e. on $\Sigma$} 
        \end{displaymath}
        for every integer $n\ge 1$. From here, we proceed in two steps.\medskip
 
        % we note that if 
        % We now define $\varphi_1$ by the graph
        % $\varphi_1 := \set{(v,w) \in L^1\times L^1 | w =
        % \varphi(v)}$.
        % By \cite[Proposition 2.17]{CoulHau2016},
        % $\varepsilon\varphi_1+(\partial\E)_{\vert_{L^{1\cap\infty}}}\varphi$
        \underline{\emph{1. Step. }} % First, consider the case $f\in D(\E\phi)$ being
        % \emph{positive}; meaning, $f\ge 0$ $\mu$-a.e. on $\Sigma$. Since
        % $\nu\phi(\cdot)+A\phi$ is $T$-accretive in $L^1$, one has
        Suppose that $u_{\lambda}\ge 0$ $\mu$-a.e. on $\Sigma$ for all
        $\lambda>0$. Then, by~\eqref{eq:pointwise-mu-conv}, Fatou's lemma applied
        to~\eqref{eq:bound-u-lambda} yields
        \begin{displaymath}
          \int_{\Sigma}j(f)\,\td\mu\le 
          \liminf_{\lambda\to0+}\int_{\Sigma}j(u_{\lambda})\,\td\mu
          \le \int_{\Sigma}j(f)\,\td\mu
        \end{displaymath}
        for every $j\in \mathcal{J}_{0}$. Choosing $j=I_{\R}$, we can
        conclude that $u_{\lambda} \to f$ in $L^{1}$ as $\lambda\to0+$
        (cf~\cite[Vol.~I., Chap.~2., Theorem~2.8.9]{MR2267655}). If
        $u_{\lambda}\le 0$ $\mu$-a.e. on $\Sigma$ for all $\lambda>0$,
        then \eqref{eq:pointwise-mu-conv} yields that $f\le 0$
        $\mu$-a.e. on $\Sigma$. Since, in particular, $-u_{\lambda}$ and
        $-f$ satisfy~\eqref{eq:bound-u-lambda}, Fatou's lemma implies
        that $-u_{\lambda} \to -f$ in $L^{1}$ as $\lambda\to0+$ and
        hence, $u_{\lambda} \to f$ in $L^{1}$ as $\lambda\to0+$.

        \underline{\emph{2. Step. }} Now, let $f\in D(\E\phi)$. By
        \eqref{eq:latice}, the positive part
        $f^{+}$ and the negative part $f^{-}$ belong to
        $D(\E\phi)$. Hence, for
        every $\lambda>0$, there are unique solutions $u_{+,\lambda}$
        and $u_{-,\lambda}$ of \eqref{eq:rangeProblemAphi} for $f$
        respectively replaced by $f^+$ and $f^{-}$. Since
        $\nu\phi(\cdot)+A\phi$ is $T$-accretive in $L^1$ and has a complete
        resolvent, for every $\lambda>0$, $u_{+,\lambda}\ge 0$ and
        $u_{-,\lambda}\le 0$ $\mu$-a.e. on $\Sigma$, and respectively
        satisfy~\eqref{eq:bound-u-lambda} with $f$ replaced by $f^+$ and
        $f^{-}$. Therefore \emph{Step 1.} implies that
        $u_{+,\lambda} \rightarrow f^{+}$ in $L^{1}$ and
        $u_{-,\lambda} \rightarrow u^{-}$ in $L^{1}$. Moreover, by the
        $T$-accretivity of $\nu\phi(\cdot)+A\phi$,
        \begin{displaymath}
          u_{-,\lambda}\le u_{\lambda}\le u_{+,\lambda}      
        \end{displaymath}    
        $\mu$-a.e.~in $\Sigma$ for every $\lambda>0$.  From this
        sandwich inequality and since $u_{+,\lambda} \rightarrow u^{+}$
        in $L^{1}$ and $-u_{-,\lambda} \rightarrow -u^{-}$ in $L^{1}$ as
        $\lambda \rightarrow 0^{+}$, one can extract from every zero
        sequence $(\lambda_{n})_{n\ge 0}$ a subsequence
        $(\lambda_{k_{n}})_{n\ge 1}$ and find a positive function
        $g\in L^{1}$ such that $\abs{u_{\lambda_{k_{n}}}}\le g$ holds
        $\mu$-a.e.~on $\Sigma$ for all $n\ge 1$. %  Thus, by Lebesgue's
        % dominated convergence theorem, it follows that
        % $u_{\lambda_{k_{n}}}\to u$ in $L^{1}$ as $n\to\infty$, and
        % thereby we have shown that the domain
        % $D((\partial \E)_{\vert_{L^{1\cap\infty}}} \varphi)$ lies dense
        % in $D(\E \varphi)\cap L^{1\cap\infty}$ with respect to the
        % $L^1$-norm topology.
\end{proof}

%\begin{proof}[Proof of Corollary~\ref{cor:densityResult}]
%    For this we refer to the subdifferential and complete
% accretivity setting of \cite{MR1164641}, noting that 
% $(\partial\E)_{\vert_{L^{1\cap\infty}}} \subseteq
%\partial_{L^{1\cap\infty}}\E$. 
%    Applying \cite[Lemma 7.1]{MR1164641} to $\partial_{
%L^{1\cap\infty}}\E$, we have that $(\partial\E)_{\vert_{
%L^{1\cap\infty}}}$ is completely accretive.    
%    Since $\E(0) = 0$ and $\E(u) \ge 0$ for all $u \in L^2$,
% it follows from the definition that $(0,0) \in \partial \E$.
% We note that since $\E$ is lower semicontinuous, again by the
% definition of the subdifferential, $\overline{\partial
% \E}^{\mbox{}_{L^2}} = \partial\E$.
%    So applying \cite[Theorem 7.4]{MR1164641}, we see that
% $\partial\E$ is $m$-completely accretive in $L^2$.
% Hence we satisfy the conditions of
% Theorem~\ref{thm:L1densityofE}.
% \end{proof}

\section{Applications}
\label{sec:applications}

We now apply our density result (Theorem~\ref{thm:L1densityofE}) to a
well-known class of doubly nonlinear operators $A\phi$ for demonstrating the role of
the hypotheses. %We begin with the following class of \emph{local} operators.

% \subsection{\texorpdfstring{Local doubly nonlinear operators}{}} In this
% first subsection
Here, we focus on doubly nonlinear operators $A\phi$ of the form
\begin{displaymath}
  A\phi(u)=-\divergence\big(\vec{a}(x,\nabla\phi(u))\big)
\end{displaymath}
equipped with \emph{homogeneous Dirichlet boundary conditions}.\medskip 

Throughout this subsection, let $\Omega\subseteq \R^d$ be an open
subset, and $\vec{a} :
 \Omega \times \R^{d}\to \R^{d}$ and $\mathcal{A}: \Omega\times
 \R^{d}\to \R$ be two Carath\'eodory functions satisfying
 \begin{align}
   \label{eq:zero-cond}
   \vec{a}(x,0)&=0 &&\\
   %&&\quad\text{for a.e. $x\in \Omega$,}\\
   \label{eq:growth-cond}
   \abs{\vec{a}(x,\xi)}&\le a_{0}(x)\,\abs{\xi}^{p-1}+a_{1}(x) 
    &&\text{ for every $\xi\in \R^d$,}\\
   \label{eq:coercivness}
    \vec{a}(x,\xi)\xi &\ge \eta \abs{\xi}^{p}-a_{2}(x) 
    &&\text{ for every $\xi\in \R^d$,}\\
   \label{eq:monotonicity-of-a}
   (\vec{a}(x,\xi_{1})&-\vec{a}(x,\xi_{1}))(\xi_{1}-\xi_{2})\ge 0 
    &&\text{ for every $\xi_1$, $\xi_2\in \R^d$,}\\
    \label{gradient:cond}
   \nabla_{\!\xi}\mathcal{A}(\cdot,\xi)&=\vec{a}(\cdot,\xi) &&
   \text{ for every $\xi\in \R^d$,}                                                               
 \end{align}
 and for a.e. $x\in \Omega$, where
 %with $\xi_{1}\neq \xi_{2}$, 
 $1<p<\infty$, $p'=p/(p-1)$, $a_{0}\in L^{\infty}$,
 $a_1\in L^{p'}$, $a_2\in L^p$, and some constant $\eta>0$.\medskip

 Further, let $W^{1,(2,p)}_{0}(\Omega)$ be the completion of
 $C^{\infty}_{c}(\Omega)$ with respect to
 $\norm{u}_{W^{1,(2,p)}_{0}}:=\norm{u}_2+\norm{\abs{\nabla u}}_p$, $u\in C^{\infty}_{c}(\Omega)
 $. It easily verifies that
 $W^{1,(2,p)}_{0}(\Omega)$ is isometric isomorphic to a closed subspace of $L^2\times
 \prod_{i=1}^{d}L^{p}$. Thus, the mixed Sobolev space
 $W^{1,(2,p)}_{0}(\Omega)$ is, in fact, reflexive and separable.\medskip
 
Under the above given hypothesis, the functional $\E : L^{2}\to
(-\infty,+\infty]$ given by
\begin{displaymath}
  \E(u):=
  \begin{cases}
    \displaystyle\int_{\Omega}\mathcal{A}(\cdot,\nabla u)\,\dx & \text{ if $u\in
      W^{1,(2,p)}_{0}(\Omega)$,}\\[5pt]
    +\infty &\text{otherwise,}
  \end{cases}
\end{displaymath}
for every $u\in L^2$ is proper, convex, and lower semicontinuous. Indeed,
$\E$ is proper since by~\eqref{eq:zero-cond}
and~\eqref{gradient:cond}, $u=0$ belongs to $D(\E)$, $\E$ is convex by
\eqref{eq:monotonicity-of-a} and \eqref{gradient:cond}, and $\E$ is lower
semicontinuous by \eqref{eq:coercivness} and by the relative weak
compactness of the unit ball of $W^{1,(2,p)}_{0}(\Omega)$. Further, the
subgradient $\partial\E$ is a well-defined mapping $\partial\E :
D(\partial\E)\to L^2$ with effective domain
\begin{align*}
  D(\partial\E)&:=
  \Big\{ u\in W^{1,(2,p)}_{0}(\Omega)\, \Big\vert\,\exists\,
  v\in L^2\text{ satisfying
                 \eqref{eq:sub-diff-characterization}}\Big\}\text{,
                 and}\\
\partial\E(u)&=v\quad\text{ for every }u\in D(\partial\E),
\end{align*}
where
\begin{equation}
  \label{eq:sub-diff-characterization}
  \int_{\Omega}\vec{a}(x,\nabla
  u)\nabla\xi\,\dx=\int_{\Omega}v\,\xi\,\dx
\end{equation}
for every $\xi\in W^{1,(2,p)}_{0}(\Omega)$. Then, the subgradient
$\partial\E$ provides a \emph{realization} in $L^2$ of the quasi-linear second-order
differential operator
\begin{displaymath}
  \partial\E(u)=-\divergence\big(\vec{a}(x,\nabla u)\big)
\end{displaymath}
equipped with \emph{homogeneous Dirichlet boundary conditions} for every
$u\in D(\partial\E)$. It follows from~\eqref{eq:zero-cond} and
\eqref{eq:sub-diff-characterization} that $(0,0)\in \partial\E$, proving
that \eqref{eq:Hyp-A0} holds. To see that $\E$ satisfies
\eqref{eq:complete-accretive-phi} and therefore, $\partial\E$ is
$m$-completely accretive in $L^2$, let $u$, $\hat{u}\in D(\E)$ and
$T\in \mathcal{P}_0$. Since $0\le T'\le 1$, it follows from the
convexity of $\mathcal{A}$ that \allowdisplaybreaks
\begin{align*}
  &\E(u+T(\hat{u}-u))+\E(\hat{u}-T(\hat{u}-u))\\
  &\qquad =
    \int_{\Omega}\mathcal{A}(x,\nabla(u+T(\hat{u}-u)))\,\dx 
    + \int_{\Omega}\mathcal{A}(x,\nabla(\hat{u}-T(\hat{u}-u)))\,\dx\\
  &\qquad =
    \int_{\Omega}\mathcal{A}(x,\nabla u+T'(\hat{u}-u)\,\nabla(\hat{u}-u))\,\dx\\ 
  &\hspace{4cm}  +
    \int_{\Omega}\mathcal{A}(x,\nabla\hat{u}-T'(\hat{u}-u)\,\nabla(\hat{u}-u))\,\dx\\
  &\qquad \le
    \int_{\Omega}(1-T'(\hat{u}-u))\mathcal{A}(x,\nabla u)\,\dx
    +\int_{\Omega}T'(\hat{u}-u)\mathcal{A}(x,\nabla \hat{u})\,\dx\\
  &\hspace{2cm}  +
    \int_{\Omega}(1-T'(\hat{u}-u))\mathcal{A}(x,\nabla\hat{u})\,\dx
    +\int_{\Omega}T'(\hat{u}-u)\mathcal{A}(x,\nabla u)\,\dx\\
  &\qquad = \E(u)+\E(\hat{u}), 
\end{align*}
proving \eqref{eq:complete-accretive-phi}. 

Next, let $\phi : \R\to \R$ be continuous, strictly monotone increasing,
surjective, and satisfy $\phi(0)=0$. Then we show that
$A=(\partial\E)_{L^{1\cap\infty}}$ and $\phi$ verify the hypotheses
\eqref{eq:sum-rule} and \eqref{eq:latice}. For $\lambda>0$, let
$[\phi^{-1}]_{\lambda}$ be the the Yosida operator of the inverse graph
$\phi^{-1}$. Since $[\phi^{-1}]_{\lambda}$ is Lipschitz continuous with
derivative $[\phi^{-1}]'_{\lambda}\ge 0$, it follows from
\eqref{eq:sub-diff-characterization}, \eqref{eq:monotonicity-of-a}, and
\eqref{eq:zero-cond} that
\begin{displaymath}
  \int_{\Omega}v\,[\phi^{-1}]_{\lambda}(u)\,\dx 
  =  \int_{\Omega}\vec{a}(x,\nabla u)\nabla
    u\,[\phi^{-1}]'_{\lambda}(u)\,\dx\ge 0
\end{displaymath}
for every $(u,v)\in A$, proving~\eqref{eq:H3-Yosida-1}. To see that $A$ and
$\phi$, in particular, satisfies~\eqref{eq:H3-Yosida-2}, consider the
sequence $(\gamma_{\varepsilon})_{\varepsilon > 0}$ of functions
$\gamma_{\varepsilon}\in C^{0,1}(\R)$ given by
    \begin{equation}
      \label{eq:gamma-varepsilon}
        \gamma_{\varepsilon}(r) = \begin{cases}
            1& \text{if $r > \varepsilon$,}\\
            \frac{r}{\varepsilon}& \text{if $-\varepsilon \le r \le \varepsilon$,}\\
            -1& \text{if $r < -\varepsilon$,}
        \end{cases}
    \end{equation}
for every $r\in \R$. Then, similarly to the previous estimate, one sees that
\begin{displaymath}
  \int_{\Omega}v\,\gamma_{\varepsilon}\left([\phi^{-1}]_{\lambda}(u)\right)\,\dx 
  =  \int_{\Omega}\vec{a}(x,\nabla u)\nabla
    u\,\gamma'_{\varepsilon}\left([\phi^{-1}]_{\lambda}(u)\right) [\phi^{-1}]'_{\lambda}(u)\,\dx\ge 0
\end{displaymath}
for every $(u,v)\in A$, $\lambda>0$, and $\varepsilon>0$. Since for
every $(u,v)\in A$, one has that $v\in L^{1\cap\infty}$ and since
\begin{displaymath}
  \lim_{\varepsilon\to
    0+}\gamma_{\varepsilon}\left([\phi^{-1}]_{\lambda}(u(x))\right)
  =\sign_0\left([\phi^{-1}]_{\lambda}(u(x))\right)
\end{displaymath}
for a.e. $x\in \Omega$, if follows from Lebesgue's dominated convergence
theorem that
\begin{displaymath}
  \int_{\Omega}v\,\sign_{0}\left([\phi^{-1}]_{\lambda}(u)\right)\,\dx\ge 0, 
\end{displaymath}
showing that \eqref{eq:H3-Yosida-2} holds and so, that
\eqref{eq:sum-rule} is satisfied. Finally, it remains to verify
\eqref{eq:latice}. First, we note that for every $u\in D(\E\phi)$, one
has that $u^+$ and $u^{-}\in D(\E\phi)$ if and only if $\phi(u^+)$ and
$\phi(u^-)\in W^{1,(2,p)}_0(\Omega)$, and by the monotonicity of $\phi$
and since $\phi(0)=0$, one has that $\phi(u^+)=[\phi(u)]^{+}$ and
$\phi(u^-)=[\phi(u)]^{-}$. Thus and since $W^{1,(2,p)}_0(\Omega)$ admits
the lattice property, it follows that $[\phi(u)]^{+}$ and
$[\phi(u)]^{-}\in W^{1,(2,p)}_0(\Omega)$ for every $u\in D(\E\phi)$,
establishing \eqref{eq:latice}.\medskip

Now, if, in addition, $\phi$ is locally Lipschitz continuous (i.e.,
$\phi\in C^{0,1}_{\loc}(\R)$), then the
set $C_c^1(\Omega)$ is a subset of $D(\E\phi)$, which is dense in
$L^1$. Thus, under these additional regularity assumption on $\phi$, it
follows from Corollary~\ref{cor:densityResult} that the following
density result holds.

\begin{theorem}\label{thm:appl}
  Let $\vec{a}: \Omega\times\R^d\to \R^d$ satisfy the
  hypotheses~\eqref{eq:zero-cond}--\eqref{gradient:cond} and
  $\phi\in C^{0,1}_{\loc}(\R)$ be strictly monotone increasing, surjective and
  satisfy $\phi(0)=0$. Further, denote by $A$ the part $\partial\E \cap L^{1\cap\infty}\times
  L^{1\cap\infty}$  in $L^{1\cap\infty}$ of the subgradient
  \begin{displaymath}
    \partial\E(u)=-\divergence(\vec{a}(\cdot,\nabla u))
  \end{displaymath}
  in $L^2$ equipped with homogeneous Dirichlet boundary conditions. Then, the
  domain $D(A\phi)$ is dense in $L^1$.
\end{theorem}

We conclude this section with the following remark.

\begin{remark}
  We emphasize that our main result (Theorem~\ref{thm:L1densityofE})
  is not restricted to be applicable merely to \emph{local} operators $A\phi$
  as in this section. In~\cite{doublyNonlinearPaper}, we applied this
  density result also to \emph{nonlocal} doubly nonlinear operators of
  the form
  \begin{displaymath}
    A\phi(u)=(-\Delta_p)^s\phi(u)
  \end{displaymath}
  for $1<p<\infty$ and $0<s<1$, $(-\Delta_p)^s$ denotes the celebrated
  fractional $p$-Laplace operator equipped with homogeneous Dirichlet
  boundary conditions. 

  Last but not least, we applied Theorem~\ref{thm:L1densityofE} to 
  differential operators $A\phi$, which were equipped with homogeneous Dirichlet boundary
  conditions. But these boundary conditions were chose only for
  simplicity. In fact, the statement of Theorem~\ref{thm:appl} remains
  true for homogeneous Neumann and Robin boundary conditions.
\end{remark}

%---------------------------------------------------------------------------
%                        References
%
%%%%%%%%%%%%%%%%%%%%%%%%%%%%%%%%%%%%%%%%%%%%%%%%%%%%%%%%%%%

%	\bibliographystyle{siam} %alpha
%	\bibliography{citations-density}

\begin{thebibliography}{10}

\bibitem{MR2033382}
{\sc F.~Andreu-Vaillo, V.~Caselles, and J.~M. Maz\'on}, {\em Parabolic
  quasilinear equations minimizing linear growth functionals}, vol.~223 of
  Progress in Mathematics, Birkh\"auser Verlag, Basel, 2004.

\bibitem{MR4041276}
{\sc W.~Arendt and D.~Hauer}, {\em Maximal {$L^2$}-regularity in nonlinear
  gradient systems and perturbations of sublinear growth}, Pure Appl. Anal., 2
  (2020), pp.~23--34.

\bibitem{MR2582280}
{\sc V.~Barbu}, {\em Nonlinear differential equations of monotone types in
  {B}anach spaces}, Springer Monographs in Mathematics, Springer, New York,
  2010.

\bibitem{BenilanCompletelyAccretiveMR1164641}
{\sc P.~B\'{e}nilan and M.~G. Crandall}, {\em Completely accretive operators},
  in Semigroup theory and evolution equations ({D}elft, 1989), vol.~135 of
  Lecture Notes in Pure and Appl. Math., Dekker, New York, 1991, pp.~41--75.

\bibitem{MR2267655}
{\sc V.~I. Bogachev}, {\em Measure theory. {V}ol. {I}, {II}}, Springer-Verlag,
  Berlin, 2007.

\bibitem{MR0348562}
{\sc H.~Brezis}, {\em Op\'erateurs maximaux monotones et semi-groupes de
  contractions dans les espaces de {H}ilbert}, North-Holland Publishing Co.,
  Amsterdam, 1973.
\newblock North-Holland Mathematics Studies, No. 5. Notas de Matem{\'a}tica
  (50).

\bibitem{MR2289546}
{\sc R.~Chill and A.~Fiorenza}, {\em Convergence and decay rate to equilibrium
  of bounded solutions of quasilinear parabolic equations}, J. Differential
  Equations, 228 (2006), pp.~611--632.

\bibitem{MR3465809}
{\sc R.~Chill, D.~Hauer, and J.~Kennedy}, {\em Nonlinear semigroups generated
  by {$j$}-elliptic functionals}, J. Math. Pures Appl. (9), 105 (2016),
  pp.~415--450.

\bibitem{doublyNonlinearPaper}
{\sc T.~Collier and D.~Hauer}, {\em A doubly nonlinear evolution problem
  involving the fractional {$p$}-{L}aplacian}.
\newblock Available on ArXiv: \href{ https://doi.org/10.48550/arXiv.2110.13401}
  {https://doi.org/10.48550/arXiv.2110.13401}, 2021.

\bibitem{CoulHau2016}
{\sc T.~Coulhon and D.~Hauer}, {\em Functional inequalities and regularizing
  properties of nonlinear semigroups - {T}heory and {A}pplications}.
\newblock Accepted for publication in the Springer Serie SMAI - Math\'ematiques
  et Applications, preprint from 2016 available
  at~\href{http://arxiv.org/abs/1604.08737}{http://arxiv.org/abs/1604.08737},
  2021.

\bibitem{MR647071}
{\sc M.~Crandall and M.~Pierre}, {\em Regularizing effects for {$u_{t}+A\varphi
  (u)=0$} in {$L^{1}$}}, J. Funct. Anal., 45 (1982), pp.~194--212.

\bibitem{MR454360}
{\sc L.~C. Evans}, {\em Differentiability of a nonlinear semigroup in
  {$L^{1}$}}, J. Math. Anal. Appl., 60 (1977), pp.~703--715.

\bibitem{MR3427974}
{\sc J.~A. Goldstein, D.~Hauer, and A.~Rhandi}, {\em Existence and nonexistence
  of positive solutions of {$p$}-{K}olmogorov equations perturbed by a {H}ardy
  potential}, Nonlinear Anal., 131 (2016), pp.~121--154.

\bibitem{MR3057171}
{\sc D.~Hauer}, {\em Convergence of bounded solutions of nonlinear parabolic
  problems on a bounded interval: the singular case}, NoDEA Nonlinear
  Differential Equations Appl., 20 (2013), pp.~1171--1190.

\bibitem{MR3369257}
\leavevmode\vrule height 2pt depth -1.6pt width 23pt, {\em The
  {$p$}-{D}irichlet-to-{N}eumann operator with applications to elliptic and
  parabolic problems}, J. Differential Equations, 259 (2015), pp.~3615--3655.

\bibitem{MR4200826}
\leavevmode\vrule height 2pt depth -1.6pt width 23pt, {\em Regularizing effect
  of homogeneous evolution equations with perturbation}, Nonlinear Anal., 206
  (2021), pp.~112245, 34.

\bibitem{Hauer2019}
{\sc D.~Hauer and J.~M. Maz{\'o}n}, {\em Regularizing effects of homogeneous
  evolution equations: the case of homogeneity order zero}, Journal of
  Evolution Equations,  (2019).

\bibitem{MR4365127}
{\sc D.~Hauer and J.~M. Maz\'{o}n}, {\em The {D}irichlet-to-{N}eumann operator
  associated with the 1-{L}aplacian and evolution problems}, Calc. Var. Partial
  Differential Equations, 61 (2022), pp.~Paper No. 37, 50.

\bibitem{alma991006084359706657}
{\sc N.~Igbida}, {\em Limite singuli\`ere de probl\`emes d'{\'e}volution
  non-lin{\'e}aires}, 1997.

\bibitem{MR3491533}
{\sc J.~M. Maz\'on, J.~D. Rossi, and J.~Toledo}, {\em Fractional
  {$p$}-{L}aplacian evolution equations}, J. Math. Pures Appl. (9), 105 (2016),
  pp.~810--844.

\end{thebibliography}

\end{document}